\documentclass[twoside,leqno]{article}

 \usepackage{amsmath}

\usepackage{geometry}                
\usepackage{graphicx}
\usepackage{graphicx}
\usepackage{amssymb}
\usepackage{epstopdf}

\usepackage{anysize}
\marginsize{.85in}{.85in}{.75in}{.75in} 

\newtheorem{thm}{Theorem}[section]

\newtheorem{defn}[thm]{Definition}
\newtheorem{example}[thm]{Example}

\newtheorem{rem}[thm]{Remark}

\newcommand{\R}{\mathbb{R}}

\newcommand{\T}{\mathcal{T}}

\newcommand{\ep}{\varepsilon}
\newcommand{\nn}{\noindent}

\newcommand{\dnrm}{\| \cdot \|}

\date{\vspace{-5ex}}

\begin{document}

\title{Existence and Uniqueness for Nonlinear Integro-Differential Equations in Real Locally Complete Spaces}
\author{Thomas E. Gilsdorf \and  Mohammad Khavanin}

\maketitle

\bigskip

\begin{abstract}  We extend existence and uniqueness results of [4] for nonlinear  integro-differential equations of Volterra type between  real locally complete  vector spaces. \\

\end{abstract}

\bigskip

\nn {AMS [2000] Subject Classification.  Primary: 45D05, 34K20;  Secondary: 46A17, 46A13.}\\

\nn {\bf Keywords}: Integro-differential equations of Volterra type, locally complete, Lyapunov - dissipative condition, strict Mackey convergence.  \\

\section{Introduction}

In [4] and [10], existence and uniqueness results are obtained for nonlinear integro-differential equations of Volterra types  of the form 

$$ x^{\prime} = H(t, x, Kx), \; \; x(0) = x_{0},  \; \; \; \; (1)$$

\nn where $H$ takes on values in a Banach space over $\R$, and $Kx$ is an integral operator depending on a continuous function $K$ having values in the same Banach space. Recently, such as in [1] and [8], generalizations of certain nonlinear problems have been extended to more general locally convex vector spaces or algebras.  In this paper we extend the main results of [4]  to the case in which the values are in a generalization of Banach spaces, specifically,  locally complete vector spaces.  Appropriate definitions are given below.\\

Throughout this paper, we consider a locally convex vector space over the field $\R$ of real numbers.  We will denote such spaces by $(E, \T)$, where $\T$ denotes the topology.  Basic properties of locally convex spaces can be found in [3], [6], and [9].  Locally convex spaces are generalizations of normed spaces, and our interest here is a class of locally convex spaces that are generalizations of Banach spaces, as described next.\\

\subsection{Locally complete spaces.}

The definition of a locally complete space relies on some information about certain kinds of bounded sets. Also, we need a way to construct linear subspaces that are normed spaces.\\

\begin{defn} In $(E, \T)$ a set $A$ is:

\begin{itemize}

\item {\bf Bounded} if, given any neighborhood $U$ of the origin, there exists a positive number $a = a_{U}$ such that $A\subset a\cdot U = \{a\cdot x : x\in U\}$. \\

\item A {\bf disk} if $A$ is both convex and balanced; i.e., if 
$$(\forall x, y \in A)(\forall s, t \in \R) \ni |s| + |t| \leq 1, \; \; sx + ty \in A.$$
\end{itemize}

 \end{defn}

\nn The unit ball of any normed space represents a set that is a bounded disk. 

\begin{defn} Let $B$ be a bounded disk in $(E, \T)$.  Denote by $E_{B}$, the linear span of $B$.  We equip $E_{B}$ with the normed topology given by the Minkowski sublinear functional of $B$ (see [9, p. 161]), defined by:

$$ (\forall x \in E_{B}) \; \; \;  \| x\|_{B} = \inf\{t \geq 0 : x \in t\cdot B  \} .$$
\end{defn}

\nn Given an arbitrary bounded set in a locally convex space, we can form the intersection of all disks that are closed and bounded, and which contain the set.  It turns out that this type of construction always leads to a set that remains bounded, c.f. [6; 7.3.4, p. 135].  Hence, without loss of generality, we may assume that a bounded set is a closed, bounded disk.

\nn The next definition represents the tool we need in order to work with normed spaces within (possibly not even metrizable) locally convex spaces.  Our definition comes from [7; 5.1.29, p. 158].\\

\begin{defn} A locally convex space $(E, \T)$ satisfies the  {\bf strict Mackey convergence condition} if for every bounded set $B$, there exists a closed, bounded disk $D$ such that $B\subset D$, and the topology of normed space $(E_{D}, \| \cdot \|_{D})$ is equivalent to the topology $\T$ on $B$.\\
\end{defn}

\nn Most of the typical spaces that occur in applications satisfy the strict Mackey convergence condition, including all metrizable locally convex spaces, countable products of such spaces, and certain inductive limits such as the space $\mathcal{D} = (\mathcal{D}, \mathcal{T}_{\mathcal{D}})$,  the space of test functions from distribution theory. For details and more information about the strict Mackey convergence condition, see [7; Section 5.1, p. 158 - 159].   Finally, we define locally complete spaces next.\\

\begin{defn}  A locally convex space $(E, \T)$ is {\bf locally complete} if for every closed bounded disk $B$, the normed space $(E_{B}, \| \cdot \|_{B})$ is complete; i.e., a Banach space. \end{defn}


\nn   In [5; 2.14, p. 20] locally complete spaces are also given the name  $c^{\infty}$, or {\it convenient} spaces.  In several references, such as [2, 7.1, p. 275], the definition of a $c^{\infty} - $ space additionally requires the space to be bornological (i.e., any linear map from $E$ to an arbitrary locally convex space $F$ is continuous if and only if it maps bounded sets to bounded sets).  We will use the definition of $c^{\infty} - $ spaces from [5; 2.14, p. 20]; i.e., locally complete spaces.   It should be noted that the structures of  $c^{\infty} - $ (locally complete) spaces have  become important in recent years due to the  use of such spaces in nonlinear distribution theory.  More details about applications of these spaces can be found in  [5],  and the references therein. \\


\nn Detailed information about locally complete spaces can be found in [7; Chapter 5].  For our purposes, the main facts are that the collection of locally complete spaces properly contains Banach spaces, and that every complete locally convex space is also locally complete, c.f. [7; Chapter 5] .  The following example illustrates a space that is strictly more general than a Banach space,  and we will apply our results to this example at the end of this paper.\\

\nn \begin{example}  Let  $\mathcal{D} = (\mathcal{D}, \mathcal{T}_{\mathcal{D}})$ denote the space of test functions from distribution theory.  \end{example}

\nn The following facts hold for $\mathcal{D}$, with relevant references given in each case:

\begin{itemize}
\item[(D1)] The space $\mathcal{D}$ can be expressed as an increasing countable union of complete metrizable locally convex spaces $(E_{n}, d_{n})$, where $d_{n}$ represents the topology of  the metric.  Moreover,  the topology $\mathcal{T} = \mathcal{T}_{\mathcal{D}}$ is the so - called inductive limit topology.  See [6 ; 12.1.1, p. 289].\\


\item[(D2)] Given any  bounded subset $A$ of $\mathcal{D}$, there exists a closed, bounded disk $D$ such that $A \subset D$, and the topology $\mathcal{T}_{\mathcal{D}}$ coincides with the normed topology $\dnrm_{D}$ on $A$. See [6; 12.1.4, p. 290]. Thus, $\mathcal{D}$ satisfies the strict Mackey convergence condition.\\
\item[(D3)] $\mathcal{D}$ is complete.  See [3; Ex 6, p. 165].\\
\item[(D4)]  $\mathcal{D}$ is not metrizable.  See [6; 12.1.5, p. 291].\\

\end{itemize}

\nn Thus, $\mathcal{D}$ is a locally complete  space that is not a Banach space.  \\

\nn  We can now state the complete problem for this paper, and the main results.\\

\section{An existence and uniqueness theorem for functions between locally complete spaces.}

We consider equation (1) given in the introduction with assumptions similar to those in [4], where the values are taken in a locally complete space.  For completeness, we state the problem here.  For the notation used below, $C[A, B]$ denotes the space of continuous functions from a set $A$ to a set $B$, and $B(a, r)$ denotes the open ball centered at $a$ of radius $r$ in a normed space.  Finally, for a scalar $\alpha$, the notation $\alpha \cdot B(a, 1)$ is equivalent to the statement: $\{ x : \|x - a\| < \alpha \}$. \\

\nn {\bf Problem 1.}\\

\nn Consider the first order nonlinear integro-differential equation of Volterra type 

\begin{equation}\label{main}  x^{\prime} = H(t, x, Kx), \; \; x(0) = x_{0} . \end{equation}

\nn Here, $x:\R \rightarrow E$, where $E = (E, \T)$ is a locally complete space over $\R$, and $Kx$ is the operator defined by 

\begin{equation}\label{eqn:Kx} (Kx)(t) = \int_{0}^{T} K(t, s, x(s)) ds, \end{equation}

with $K$ and $H$ satisfying the following:

 $$ K\in C[\R^{2}\times E, E],  \; \; H\in C[\R \times E \times E, E].$$
 
\nn The differentiation in $E = (E, \T)$ is defined as done for general locally convex spaces, such as in [11], or  [2; 10.2, p. 279].\\

\nn The following is inspired by Theorem 3.1 of [4].   Our result here is for functions from a locally complete space  to the same locally complete space that satisfy some boundedness conditions outlined below.  For consistency, we have chosen to use notation that coincides as closely as possible to that of [4]. \\


\nn {\bf Theorem 1.}  Assume there exists a closed, bounded disk $B\subset E$ such that for $J = [0, T]$, and some $K_0 > 0, H_{0} > 0$, 

\begin{equation}\label{eqn:A1} K\left( J^{2}\times B \right) \subset K_{0}\cdot B ; \tag{A1} \end{equation}

\begin{equation}\label{A2} H\left(J\times B \times  K_{0}\cdot B \right) \subset H_{0}\cdot B ,\tag{A2} \end{equation}

\vspace{0.15in}

\nn $K$ is Lipschitz in the third argument with respect to the norm $\dnrm_{B}$ on  $E_{B}$; in particular, 



 \begin{equation}\label{A3} (\exists k_{1} > 0) \ni \; \|K(t, s, u) - K(t, s, \overline{u}) \|_{B} \leq k_{1}\|u - \overline{u}\|_{B} , \tag{A3} \end{equation}
 
 on $J^{2} \times B$ ,\\
 
\nn $H$ is locally Lipschitz; that is, for any $(t, x, y) \in J\times B \times K_{0}\cdot B$, there exist $\delta = \delta(t,x,y) > 0$,  $L = L(t, x, y) > 0$,  and neighborhoods $U_{x}$ of $x$ and $U_{y}$ of $y$ within $J\times B \times K_{0}\cdot B$ such that in $(E_{B}, \| \cdot \|_{B})$,
 
\begin{equation}\label{A4}  \| H(t, x_{1}, y_{1}) - H(t, x_{2}, y_{2})\|_{B} \leq L\left( \|x_{1} - x_{2}\|_{B} + \|y_{1} - y_{2}\|_{B} \right) , \tag{A4} \end{equation}

\nn for $(t, x_{1}, y_{1}), (t, x_{2}, y_{2})\in J\times U_{x}\times U_{y}$.\\

\nn {\it Then}: 

\nn {\bf (a)}: There exists $\eta > 0$ such that equation ~(\ref{main}) has a unique solution on $J_{0} = [0, \eta ]$ in $(E_{B}, \| \cdot \|_{B})$.\\

\nn {\bf (b)}:  The sequence of approximations that converge to the unique solution in $(E_{B}, \| \cdot \|_{B})$ from part (a) converges to the solution with respect to the topology $\T$ of $E$.\\

\bigskip


\nn {\it Proof}:  We start by rewriting some of the assumptions in terms of norms.  By assumption (A1), $K\left( J^{2}\times B \right) \subset K_{0}\cdot B$ implies that on $J^{2}\times B$, 

$$ \| K(t, s, x)\|_{b} \leq K_{0}. \; \; \; \; (A1) $$

  The assumption (A2) implies that on $B^{\circ} $, that is, on $\{ x \in E_{B} : \|x \|_{B} < 1\}$, and on $K_{0}\cdot B = \{ x \in E_{B} : \|x \|_{B} \leq K_{0}\}$ we have 
  
$$\|H(t, x, y)\|_{B} \leq H_{0}.  \; \; \; \; (A2) $$  

\bigskip

\nn We will prove that $\eta = \min\left\{ T, \frac{1}{2H_{0}} \right\}$ is the desired value.

\nn  For $t = 0$, $x = x_{0}$, $y = 0$, choose $\sigma_{1}, \; \gamma_{1}, \; \overline{\gamma_{1}} > 0$ such that $\sigma_{1}H_{0} <  \overline{\gamma_{1}} \;$, $\sigma_{1}K_{0} <  \gamma_{1}$, and for which the Lipschitz inequality of  ~(\ref{A4}) holds on $R_{1} = Rx_{0}$, where

$$R_{1} =  Rx_{0} = [0, \sigma_{1}] \times B\left(x_{0},  \gamma_{1} \right) \times B\left(0, \overline{\gamma_{1}} \right) .$$

\nn By known methods such as Schauder's fixed point theorem or successive approximations, a unique solution $x(t)$ of (1)  can be found on $[0, \sigma_{1} ]$ as a limit of a sequence $(x_{m}) = (x_{m}(t))$, with respect to the norm $\dnrm_{B}$.  We now enlarge the interval of solution as follows. Let

$$ x_{\sigma_{1}} = x_{0} +  \int_{0}^{\sigma_{1}} H(s, x(s), (Kx)(s)) ds .$$

\nn For $t = \sigma_{1}$ and $x = x_{\sigma_{1}}$, let

$$y_{\sigma_{1}} =  \int_{0}^{\sigma_{1}} K(t, s, x(s)) ds. $$

\nn There exists $R_{2} = R_{x_{\sigma_{1}}}$ for which the Lipschitz inequality of (A4) holds, given by 

$$ R_{2} =  [\sigma_{1}, \sigma_{2}] \times B\left( x_{\sigma_{1}},  \gamma_{2} \right) \times B \left( y_{\sigma_{1}}, \gamma_{2} \right) ,$$

\nn where 

$$ \left( \sigma_{2} - \sigma_{1}\right) H_{0} < \gamma_{2}  \left( \sigma_{2} - \sigma_{1}\right)K_{0} < \overline{\gamma_{2}} ,$$

\bigskip

\nn and such that the Lipschitz inequality of  ~(\ref{A4}) holds on $R_{2}$.  It follows that we can prove the existence of a unique solution $x(t)$ on $[0, \sigma_{1} + \sigma_{2}]$.  We again denote the sequence of successive approximations by $(x_{m})$.  Let $S$ be the set of all unique solutions $x(t)$ to (1) on an interval $[0, \alpha]$ for $\alpha \leq T$.  It is easy to prove that a partial ordering of $S$ is given by set inclusion of intervals, and that Zorn's Lemma applies.  Thus, we conclude that there is a maximal element, that is, there exists a unique solution to (1), on $[0, \eta ]$. This proves part (a). To prove part (b), by [7; 3.2.2, p. 82], the topology of the norm $\dnrm_{B}$ is stronger than the topology $\T$ on the vector space $E_{B}$.  We conclude that if $(x_{m})$ converges to the unique solution on $[0, \eta ]$ with respect to $\dnrm_{B}$, then $(x_{n})$ converges to the unique solution $x = x(t)$ on $[0, \eta ]$ with respect to $\T$ as well.   $\; \; \Box$\\

\bigskip

\nn {\bf Remark.}  We proved this result by using the unit ball $B$ of $(E_{B}, \dnrm_{B})$.  In general, one can prove Theorem 1 on a ball of radius $N$ in $(E_{B}, \dnrm_{B})$.  In this case, the details follow as in the proof of Thm 3.1 of [4], with $\eta = \min\left\{ T, \frac{N}{2H_{0}} \right\}$.\\


\section{Existence and uniqueness under a Lyapunov - dissipative condition}

\nn The result that follows generalizes Theorem 3.2 of [4] to locally complete spaces, under the assumptions of our previous theorem.  \\

\nn Within the context of a locally complete space $E$ of Theorem 1, we say that $H(t, x, y)$ satisfies a {\bf Lyapunov - dissipative condition} if items (i) - (iii) below are satisfied:\\

\begin{itemize}

\item[(i)]  \begin{equation}\label{A5}  \; \; \;  V \in C\left[ J \times B \times B, \R^{+} \right], \; \; V(t, x, x) \equiv 0, \;  \; \; V(t, x, y) > 0 ,\tag{A5} \end{equation}
if $x \neq y$ , for 

$$ (t, x), (t, y) \in J \times B^{\circ}, \; \;  K_{0} \cdot B \subset B \subset E_{B} ,$$

\nn with $L > 0$ such that

$$ \left| V(t, x, y) - V(t, x_{1}, y_{1}) \right| \leq L\cdot \left( \|x - x_{1}\| + \|y - y_{1} \|_{B} \right) .$$

\item[$(i)^{\prime}$] If $(x_{n})$ and $(y_{m})$ are sequences in $B$ such that $lim_{m, n \rightarrow \infty} V(t, x_{n}, y_{m}) = 0$, then $lim_{m, n \rightarrow \infty} ( x_{n} - y_{m}) = 0$ in the topology of $E$.\\

\item[(ii)] The following derivative relation holds:

$$D\left( V(t, x, y) \right) = \lim_{h \rightarrow 0^{+}} \frac{1}{h} \left\{ V(t, x, y) - V(t - h, x - hH(t, x, Kx) , y - hH(t, y, Ky) \right\} ,$$

\nn and we have 

$$ D\left( V(t, x, y) \right) \leq g\left( t, \; V(t, x, y), \; \int_{0}^{t} S(t, s, V(s, x(s), y(s)) ds \; \right) , $$

for $t \in J$ and $x, y \in C[J, B]$, with $S\in C\left[ J \times J \times \R^{+}, \R \right]$, $\left| S(t, s, V)\right| \leq S_{0}$ on $J \times J \times \R^{+}$; moreover, $S$ satisfies condition (L4) of [4].\\

\item[(iii)]  The function $g$ satisfies: $g \in U^{\ast}$ of [4] with respect to $S$ and $t_{0} = 0$.\\
\end{itemize}

\bigskip

\begin{thm} Assume the hypotheses of Theorem 1, in particular, assumptions (A1) - (A3).  Further, assume the  Lyapunov dissipative condition (A5) and that the space $(E, \mathcal{T})$ satisfies the strict Mackey convergence condition.  Then there exists $\eta > 0$ such that equation (1) has a unique solution on $J_{0} = [0, \eta]$. \\ 

\end{thm}

\bigskip

\nn {\it Proof:}  As in the proof of Theorem 3.2 in [4], we construct a sequence $\{ x_{n}(t)\}$ of $\ep_{n} - $ approximations on the interval $J_{0} = [0, \eta]$, where $0 < \ep_{n} < 1$ and $\ep_{n} \longrightarrow 0$ as $n \longrightarrow \infty$.  To finish the proof, it will suffice to prove that the sequence converges to a continuous function $x(t)$ in the topology $\mathcal{T}$ of $E$; the proof that $x(t)$ is the unique solution follows from the same arguments as the proof of Theorem 3.2 in [4].\\

\nn By [4; Thm 2.2, p. 94] and the assumption that $g \in U^{\ast}$, the arguments of Step II of [4; Thm 3.2, p. 101] apply to conclude that 

$$ \lim_{n \rightarrow \infty , m\rightarrow \infty} \left[ V(t, x_{n}(t), x_{m}(t) \right] = 0 ,$$

\nn in the topology $\mathcal{T}$ of $E$, for any $t \in J_{0}$.   By assumption of the strict Mackey convergence condition, there is a closed, bounded disk $D$, with $B \subset D$, and for which $\mathcal{T}_{D}$ is equivalent to the normed topology of $(E_{D}, \| \cdot \|_{D})$, on the set $B$. By local completeness, $(E_{D}, \| \cdot \|_{D})$ is a Banach space, and we may apply the arguments from Step II of [4; p. 101], to conclude that $(x_{n}(t))$ is uniformly Cauchy in $(E_{D}, \| \cdot \|_{D})$.  Hence,  $(x_{n}(t))$ converges in the space $(E_{D}, \| \cdot \|_{D})$ to a continuous function $x(t)$.  Finally, by the equivalence of the norm $\dnrm_{D}$ to the topology of $E$ on the set $B$, we conclude that $(x_{n}(t))$ converges to  $x(t)$ in the space $(E, \mathcal{T}). \; \; \Box$ \\

\begin{rem}  In view of items (D1) - D(4), the results of Theorem 3.3 hold for the space $\mathcal{D} = (\mathcal{D}, \mathcal{T}_{\mathcal{D}})$, the space of test functions from distribution theory.  \end{rem}

\bigskip

\begin{center} {\sc References} \end{center}

\bigskip

\nn [1] Bosch, C., Garc\'{i}a, A., G\'{o}mez - Wulschner, C., Hern\'{a}ndez - Linares, S. {\it Equivalents to Ekeland's variational principle in locally complete spaces}.  Sci., Math. Jpn., {\bf 72}, (2010), no. 3, 283 - 287.\\

\nn [2]  Fr\"{o}licher, A. {\it Axioms for convenient calculus}.  Cah. de Topol. et G\'{e}om. Diff\'{e}er. Cat\'{e}., {\bf 45}, (2004), no. 4, 267 - 286.\\

\nn [3]  Horv\'{a}th, J. \underline{Topological Vector Spaces and Distributions, Vol. I}.  Addison - Wesley, 1966.\\

\nn [4]  Hu, S., Wan, Z., Khavanin, M. {\it On the existence and uniqueness for nonlinear integro - differential equations}.  Jour. Math Phy. Sci., {\bf 21}, no. 2, (1987), 93 - 103.\\

\nn [5]  Kriegl, A., Michor, P. W. \underline{The Convenient Setting Global Analysis}.  AMS - Surveys and Monographs, {\bf 53}, (1997).\\

\nn [6]  Narici, L., Beckenstein, E. \underline{Topological Vector Spaces}.  M. Dekker, (1985).\\

\nn [7] P\'{e}rez Carreras, P., Bonet, J. \underline{Barrelled Locally Convex Spaces}.  North Holland, 1987.\\

\nn [8]  Stojanovi\'{c}, M. {\it System on nonlinear Volterra's integral equations with polar kernel and singularities}.  Nonlinear Anal., {\bf 66}, no.7, (2007), 1547 - 1557.\\

\nn [9]  Swartz, C.  \underline{Functional Analysis}.  M. Dekker, (1992).\\

\nn [10]  Wan, Z.  {\it Existence and uniqueness of solutions of nonlinear integro - differential equations of Volterra type in a Banach space}.  Appl. Anal., {\bf 22}, (1986), 157 - 166.\\

\nn [11]  Yamamuro, S. \underline{A theory of differentiation in locally convex spaces}. Mem. Amer. Math. Soc.,  {\bf 17} (1979), no. 212.

\vspace{0.75in}



%



\nn {\sc Department of Mathematics\\ University of North Dakota\\ 101 Cornell Street, Mail Stop 8376\\ Grand Forks, ND 58202-8376\\ USA}\\

\nn {\it E-mail}: \/   thomas.gilsdorf@und.edu,\/  mohammad.khavanin@und.edu

\end{document}